\documentclass[aos,MSNbibl,citesort,dvips]{arximspdf}
\usepackage{mathbh}
\usepackage{graphicx}

\doi{10.1214/11-AOS929}
\volume{39}
\issue{6}
\pubyear{2011}
\firstpage{2974}
\lastpage{3002}

\makeatletter

\newtheorem{corollary}{Corollary}
\newtheorem{lemma}{Lemma}

\newproclaim{assumption}{Assumption}
\newproclaim{definition}[assumption]{Definition}

\newtheorem{theorem}{Theorem}[section]
\newtheorem{prop}{Proposition}[section]

\newcommand{\ind}{\mathbh{1}}

\makeatother

\begin{document}
\begin{frontmatter}

\title{Accurate emulators for large-scale computer experiments}
\runtitle{Accurate emulators for large-scale computer experiments}

\begin{aug}
\author[A]{\fnms{Ben} \snm{Haaland}\thanksref{t1}\ead[label=e1]{benjamin.haaland@duke-nus.edu.sg}}
\and
\author[B]{\fnms{Peter Z. G.} \snm{Qian}\corref{}\thanksref{t2}\ead[label=e2]{peterq@stat.wisc.edu}}
\runauthor{B. Haaland and P. Z. G. Qian}
\affiliation{Duke-NUS Graduate Medical School and National University
of Singapore, and University of Wisconsin, Madison}
\address[A]{Centre for Quantitative Medicine\\
Office of Clinical Sciences\\
Duke-NUS Graduate Medical School\\
Singapore 169857\\
and\\
Department of Statistics and Applied Probability\\
National University of Singapore\\
Singapore 117546\\
\printead{e1}}
\address[B]{Department of Statistics\\
University of Wisconsin, Madison\\
Madison, Wisconsin 53706\\
USA\\
\printead{e2}} 
\end{aug}

\thankstext{t1}{Supported by Singapore National Medical Research Council
Grant IRG10nov006 and an NUS Initiative to Improve Health in Asia,
Global Asia Institute grant.}

\thankstext{t2}{Supported by NSF Grant CMMI 0969616, NSF Career Award
DMS-10-55214 and
an IBM Faculty Award.}

\received{\smonth{5} \syear{2011}}
\revised{\smonth{9} \syear{2011}}

%
\begin{abstract}
Large-scale computer experiments are becoming increasingly important in
science. A multi-step procedure is introduced to statisticians for
modeling such experiments, which builds an accurate interpolator in
multiple steps.
In practice, the procedure shows substantial improvements in overall
accuracy, but its theoretical properties are not well established. We
introduce the terms nominal and numeric error and decompose the overall
error of an interpolator into nominal and numeric portions. Bounds on
the numeric and nominal error are developed to show theoretically that
substantial gains in overall accuracy can be attained with the
multi-step approach.
\end{abstract}

%
\begin{keyword}[class=AMS]
\kwd[Primary ]{41A17}
\kwd[; secondary ]{65M12}
\kwd{65G50}.
\end{keyword}
\begin{keyword}
\kwd{Computer experiment}
\kwd{emulation}
\kwd{interpolation}
\kwd{Gaussian process}
\kwd{large-scale problem}
\kwd{multi-step procedure}
\kwd{numerical technique}
\kwd{radial basis function}
\kwd{reproducing kernel Hilbert space}.
\end{keyword}

\end{frontmatter}

\section{Introduction}\label{intro}

Computer experiments use complex mathematical mo\-dels implemented in
large computer codes to study real systems. In many situations, a
physical experiment is not feasible because it is unethical,
impossible, inconvenient or too expensive. A mathematical model of the
system can often be developed and input/output pairs can be produced
with the help of computers. Typically, the input/output pairs are
expensive in the sense that they require a great deal of time and
computing to obtain and they are nearly deterministic in the sense that
a particular input will produce almost the same output if given to the
computer experiment on another occasion. Computer experiments are
widely used in systems biology, engineering design, computational
biochemistry, climatology and epidemiology and their pervasiveness\vadjust{\goodbreak} in
science, engineering and medicine is only growing. When using a
computer experiment to study a real system, a thorough exploration of
the surface is typically wanted. However, obtaining input/output pairs
is often too expensive for a complete exploration. %
A solution is to evaluate the computer experiment at several
well-distributed data sites given by a \textit{space-filling design}
\mbox{\cite{mcb79,owen92,tang93,flwz00,niederreiter,owen1995,jh08}}.
Then build an interpolator which can be used as a stand-in, or
emulator, for the actual computer experiment. %
The thorough exploration of the complex surface can then be carried out
on the interpolator. Excellent overviews on data collection and
modeling for computer experiments can be found in \cite
{SSW1989,SWMW1989,CMMY1991,koehowen96,santner,fang}.

To emulate the output from a computer experiment, Gaussian process~(GP)
models or reproducing kernel Hilbert space (RKHS) interpolators are
often used. These interpolators have a simple form and control the
smoothness of the emulator.
In particular, let $f$ denote the output of a run of the computer
experiment, so that the functional link between input $x$ and output
$y$ is $y=f(x)$.
Take $\Phi\dvtx\Omega\times\Omega\to\mathbb{R}$ to be symmetric in
its two
arguments and positive definite. The kernel $\Phi$ is \textit{positive
definite} on a domain of interest~$\Omega$~if
\[
\sum_{i=1}^n\sum_{j=1}^n\alpha_i\alpha_j\Phi(x_i,x_j)>0
\]
for every nonzero $\alpha\in\mathbb{R}^n$ and distinct $\{
x_1,\ldots
,x_n\}\subseteq\Omega$. Then, given distinct input sites $X=\{
x_1,\ldots
,x_n\}$, the GP or RKHS interpolator has the simple form
\[
\mathcal{P}(x)=\sum_{i=1}^n\alpha_i\Phi(x,x_i),
\]
where $\alpha$ has $A_X \alpha=f|_X$, $A_X=\{\Phi(x_i,x_j)\}$ and
$f|_X=(f(x_1)\cdots f(x_n))'$.
Associated with each symmetric, positive definite kernel is exactly one
Hilbert space of functions whose norm, in the case that the kernel is
smooth, measures both size and smoothness. For a particular kernel
$\Phi
$, this associated function space will be called its \textit{native
space} and will be denoted $\mathcal{N}_\Phi(\Omega)$. Native spaces
will be discussed further in Section \ref{nominalsection}. The
smoothness of the emulator is controlled in the sense that the RKHS
interpolator has the smallest possible native space norm of any
function interpolating $f|_X$ \cite{fasshauer,wahba,wendland2005}.
It is worth noting that the GP models often used in practice to build
emulators for computer experiments are essentially a special case of
RKHS emulators.
In the GP context, the kernel $\Phi$ is a, possibly scaled, correlation
function.
In the case that a nonzero mean function $\hat{\mu}$ is estimated in
the GP model, the interpolator is actually the sum of this estimated
mean function and an RKHS interpolator of the residual $(f-\hat{\mu})|_X$.
Here, we consider \textit{translation invariant}, or stationary, kernels
so that $\Phi$ is a function of only the difference between its
arguments. Hereafter, $\Phi(x,y)$ will be written as $\Phi(x-y)$.
Note that the connection between Gaussian processes and RKHS was also
discussed in\vadjust{\goodbreak}
\cite{wahba}.




Many of the systems which scientists, engineers and medical researchers
use computer experiments to study exhibit extremely complex behavior in
portions of the input space. To discover and understand these regions
requires a large-scale computer experiment with many input sites which
are potentially very near one another. Unfortunately, most methods for
building emulators, including RKHS and GP interpolators, suffer from
increasingly poor predictive accuracy due to numerical problems as the
number of observations of the computer experiment becomes larger.
Throughout, we refer to large-scale computer experiments as those with
a large number of runs.
Such experiments appear frequently
in various fields such as aerospace engineering \cite{Boeing},
information technology \cite{IBM}, biology, high-energy physics,
nanotechnology and security.
The essential difficulty in emulation of a large-scale computer
experiment is that as input sites become nearer to one another the
problem of finding an interpolator becomes ill-conditioned and so less
amenable to accurate calculation. Several techniques exist for
numerically stabilizing kernel-based interpolators, including adding a
nugget effect
\cite{santner,LBH2006}, using compactly supported kernels \cite
{Gneiting2002,fasshauer},
covariance tapering \cite{kauf2008}, decomposing the correlation matrix
\cite{booker} and
approximating likelihoods \cite{stein2004}. The multi-step procedure
\cite{floater} described below also addresses the vital issue of
numerical stability and can be used alone or in concert with additional
numeric measures such as those mentioned above.

The multi-step procedure is not new to the field of applied
mathematics, yet the exposure of statisticians to this method is
relatively limited. Further, while the procedure often improves overall
predictive accuracy substantially in practice, minimal work has been
done on its theoretical properties \cite{fasshauer}. Notable exceptions
include \cite{narcowich}, \cite{fasshauer2} and \cite{HaLe02}. The existing
theoretical work in the literature examines numerical accuracy in a
relatively qualitative manner. Here, we introduce the concepts of \textit
{nominal} and \textit{numeric} accuracy. Nominal accuracy refers to the
accuracy which would be attained if computations could be performed
without floating point rounding. Numeric accuracy refers to how close
computed quantities are to their corresponding nominal counterparts.
Then, we introduce a decomposition of the error of an interpolator into
nominal and numeric portions. This gives a complete description of the
computed interpolator's error while separating the contributing sources
of error to allow for more straight-forward analysis. Bounds on the
numeric and nominal error of the multi-step interpolator are developed.
The numeric bound is the only complete, rigorous bound on the numeric
error of the multi-step interpolator. The result is very general and
makes very few assumptions about the kernels used in different steps.
The nominal bound is similar to the error bound developed in \cite
{narcowich}, but more general in that it allows the kernels at
different stages to be re-scaled in a~flexible manner. In practice, the
kernel re-scalings can have a large impact on accuracy.

%
%

\section{Multi-step interpolator}\label{multi-stepSection}

The multi-step procedure explored here is a~generalization of the
procedure introduced in
\cite{floater}.
Their idea was to form well-spread nested subsets of the data. Then
interpolate the first subset using a wide kernel and form residuals of
this interpolator on the next subset. The residuals are then
interpolated using a narrower kernel and the current stage, and
previous stage interpolators are added together, giving an interpolator
on the larger subset. This procedure is repeated an appropriate number
of times, at each stage updating the interpolator, until an
interpolator of the complete data is obtained.
We introduce a separation of the error into nominal and numeric
portions and derive bounds on each type of error.
We adopt a slightly different notation than
\cite{floater}.
Let $f$ denote the unknown function to be interpolated and $\Omega
\subseteq\mathbb{R}^d$ denote the domain of interest.
Throughout, the following assumption is made about the
kernel~$\Phi$.\vspace*{-2pt}

\begin{assumption}\label{kernelAssumption}
The kernel $\Phi$ is continuous, positive definite and translation
invariant.\vspace*{-2pt}
\end{assumption}

Note that with minor modifications, the development and
results in Sections \ref{intro}--\ref{numSec3} only require that
$\Phi$
is positive definite.

In the below description of the multi-step interpolation procedure, $J$
denotes the number of stages, and $\Phi_j$ denotes the 
kernel used for interpolation in
stage $j$.
Now, take
%
\begin{equation}\label{eqnnested-design}
X_1\subset\cdots\subset X_J=X\vspace*{-2pt}
\end{equation}
and initialize $\mathcal{P}^0\equiv0$. Then, for $j=1,\ldots,J$, let
%
\begin{eqnarray}\label{multi-step}
\mathcal{P}^j(x)&=&\sum_{u=1}^{n_j}\alpha^j_u\Phi
_{j}(x-x_u),\nonumber\\[-2pt]
\alpha^j&=&A^{-1}_{X_j,\Phi_j}\Biggl(f-\sum_{k=0}^{j-1}\mathcal
{P}^k\Biggr)\Bigg|_{X_j},\nonumber\\[-9pt]\\[-9pt]
A_{X_j,\Phi_j}&=&\{\Phi_{j}(x_u-x_v)\},\qquad u,v=1,\ldots,n_j,\nonumber\\
n_j&=&\operatorname{card} X_j.\nonumber\vspace*{-2pt}
\end{eqnarray}
Then the multi-step interpolator,
%
\begin{equation}\label{multi-step2}
\mathcal{P}(x)=\sum_{j=1}^J\mathcal{P}^j(x)\vspace*{-2pt}
\end{equation}
satisfies the interpolation conditions $\mathcal{P}(x_u)=f(x_u)$,
$u=1,\ldots,n$, where $n=\operatorname{card} X$.
Here, $X$ is the complete set of input sites.
The results in this article indicate that the best performance will be
achieved if each of the nested designs, $X_1,\ldots,X_J$, are chosen to
have well-separated data sites, uniform low-dimensional projections and
small data-free regions.
Note\vspace*{1pt} that $\alpha^j$ should not be calculated\vadjust{\goodbreak} using the formula
$A^{-1}_{X_j,\Phi_j}(f-\sum_{k=0}^{j-1}\mathcal{P}^k)|_{X_j}$, but
instead as the solution to the linear system $A_{X_j,\Phi_j}\alpha
^j=(f-\sum_{k=0}^{j-1}\mathcal{P}^k)|_{X_j}$.
In general, the solution to the linear system is subject to smaller
numeric error.
Also, in the situation where $n$ is large and $A_{X_j,\Phi_j}$ is
sparse due to memory constraints, $A_{X_j,\Phi_j}^{-1}$ will often be
too dense to be stored.

It is commonly the situation that each kernel $\Phi_j$ depends on
parameters $\Theta_j$.
For example, in Section \ref{nominalsection} it is assumed that~$\Phi
_j$ is a known kernel~$\Psi_j$ whose inputs $x-y$ are re-scaled by a
matrix $\Theta_j$, so that $\Phi_j(x-y)=\Psi_j(\Theta_j(x-y))$.
The form of the underlying kernels $\Phi_j$ is often fixed in advance
to achieve an interpolator with prespecified smoothness and numerical
properties.
In particular, the results in Sections \ref{secnum} and \ref
{nominalsection} indicate that smoother underlying kernels have better
nominal properties and worse numeric properties, as defined in (\ref
{errorSeparation}), and vice versa.
The accuracy of the interpolator can depend significantly on the choice
of parameter values.
A few possible criteria for choosing the parameters $\Theta_j$ are
cross-validation, maximum likelihood and sparsity of the interpolation matrices.
Most procedures for choosing the $\Theta_j$ are simplified by
considering each stage sequentially.
In particular, $\Theta_j$ can be chosen to minimize the
cross-validation error, maximize the likelihood or restrain the number
of nonzero entries in the interpolation matrix $A_{X_j,\Phi_j}$ at
stage $j$.
For smaller problems, where a dense $A^{-1}_{X_j,\Phi_j}$ can be
stored, the short-cut formula in (\ref{shortcut}) can be used to make
leave-one-out cross-validation computationally efficient.
For larger problems, an option such as 10-fold cross-validation is more
appropriate.
If the residuals from the previous stage $(f-\sum_{k=0}^{j-1}\mathcal
{P}^k)|_{X_j}$ are modeled as a GP, then maximum likelihood can be used
to choose the parameters~$\Theta_j$.
Maximizing the likelihood at each stage is equivalent to minimizing
%
\begin{equation}\label{mlcriterion}
n_j\operatorname{log} \Biggl[\frac{1}{n_j}\Biggl(f-\sum_{k=0}^{j-1}\mathcal
{P}^k\Biggr)'\Bigg|_{X_j} \alpha^j\Biggr]+ \operatorname{log} \operatorname{det}(A_{X_j,\Phi_j}).
\end{equation}
Restricted maximum likelihood estimates can be obtained by replacing
the~$n_j$ in the objective function (\ref{mlcriterion}) by
$n_j-n_{j-1}$, with $n_0=0$.
For large problems, a storage and computation efficient algorithm such
as \cite{barry} should be used in calculating $\operatorname{log}
\operatorname{det}(A_{X_j,\Phi_j})$.
For very large problems, memory constraints demand that the sparsity of
$A_{X_j,\Phi_j}$ be considered.
One possibility for compactly supported kernels is to choose \textit
{fixed} $\Theta_j$ to ensure that the number of nonzero entries in
$A_{X_j,\Phi_j}$ is manageable as in (\ref{fixedTheta}).
Another possibility is to incorporate a penalty for nonsparsity into
the objective function such as~(\ref{mlcriterion}).


If the error at stage $j$, $f-\sum_{k=0}^{j-1}\mathcal{P}^k$, is
modeled as a GP, then confidence intervals on the function's values
$f(x)$ can be obtained in much the same manner as a single stage
interpolator \cite{wahba1983}.
In particular,
model the output~as
\[
f(x)=\sum_{j=1}^J Z_j(x),\vadjust{\goodbreak}
\]
where the $Z_j$ are mean zero Gaussian processes with
$\operatorname{Cov}(Z_j(x_1),Z_j(x_2))=\sigma^2_j\Phi_j(x_1-x_2)$.
Note that the $Z_j$ are \textit{not} independent.
For point sets $X$ and~$Y$, denote the $\operatorname{card} X\times\operatorname{card}
Y$ matrix of pairwise kernel evaluations of points in $X$ and $Y$ as
%
\begin{equation}\label{notation}
\Phi(X-Y)=\{\Phi(x_u-y_v)\},
\end{equation}
where $x_u\in X$, $y_v\in Y$.
Take $Z_0\equiv0$ to simplify the development below.
Conditional on $f|_{X_J},Z_1,\ldots,Z_{J-1}$,
%
\begin{eqnarray}\label{cond1}\quad
&&f(x)-\sum_{j=0}^{J-1}Z_j(x)
\sim\mathcal{N}\Biggl(\Phi_J(X_J-x)'A_{X_J,\Phi_J}^{-1}\Biggl(f-\sum
_{j=0}^{J-1}Z_j\Biggr)\Bigg|_{X_J},\nonumber\\[-3pt]
&&\hspace*{111pt}\sigma^2_J\bigl(\Phi_J(0)-\Phi_J(X_J-x)'A_{X_J,\Phi
_J}^{-1}\Phi(X_J-x)\bigr)\Biggr)\nonumber\\[-10pt]\\[-10pt]
&&\quad\Longrightarrow\quad f(x)
\sim\mathcal{N}\Biggl(\Phi_J(X_J-x)'A_{X_J,\Phi_J}^{-1}\Biggl(f-\sum
_{j=0}^{J-1}Z_j\Biggr)\Bigg|_{X_J} +\sum_{j=0}^{J-1}Z_j(x),\nonumber\\[-3pt]
&&\hspace*{125.5pt}\sigma^2_J\bigl(\Phi_J(0)-\Phi_J(X_J-x)'A_{X_J,\Phi
_J}^{-1}\Phi(X_J-x)\bigr)\Biggr).
\nonumber
\end{eqnarray}
Let $\tilde{X}_J=\{X_J, x\}$.
Then, conditional on $f|_{X_J},Z_1,\ldots,Z_{j-1}$
%
\begin{eqnarray}\label{cond2}\quad
%
Z_j|_{\tilde{X}_J}&\sim&\mathcal{N}\Biggl(\Phi_j(X_j-\tilde
{X}_J)'A_{X_j,\Phi
_j}^{-1}\Biggl(f-\sum_{k=0}^{j-1}Z_k\Biggr)\Bigg|_{X_j},\nonumber\\[-10pt]\\[-10pt]
&&\hphantom{\mathcal{N}\Biggl(}\sigma^2_j\bigl(\Phi_j(\tilde{X}_J-\tilde{X}_J)-\Phi
_j(X_j-\tilde{X}_J)'A_{X_j,\Phi_j}^{-1}\Phi_j(X_j-\tilde
{X}_J)\bigr)\Biggr).\nonumber
\end{eqnarray}
Note that the distribution in (\ref{cond2}) is singular and $\Phi
_j(\tilde{X}_J-\tilde{X}_J)=A_{\tilde{X}_J,\Phi_j}$ in the notation of
(\ref{multi-step}).
The first $n_j$ components of these conditional distributions are
trivial and given by
\[
Z_j|_{X_j}=\Biggl(f-\sum_{k=0}^{j-1}Z_k\Biggr)\Bigg|_{X_j},
\qquad j=1,\ldots,J.
\]
The remaining $n_J-n_j+1$ components have the nontrivial distribution,
conditional on $f|_{X_J},Z_1,\ldots,Z_{j-1}$, given by
\begin{eqnarray*}
&&Z_j|_{\tilde{X}_J\setminus X_j}\sim\mathcal{N}\Biggl(\Phi_j(X_j-\tilde
{X}_J\setminus X_j)'A_{X_j,\Phi_j}^{-1}\Biggl(f-\sum
_{k=0}^{j-1}Z_k\Biggr)\Bigg|_{X_j},\\[-3pt]
&&\quad\sigma^2_j\bigl(\Phi_j(\tilde{X}_J\setminus X_j\,{-}\,\tilde
{X}_J\setminus X_j)\,{-}\,\Phi_j(X_j\,{-}\,\tilde{X}_J\setminus X_j)'A_{X_j,\Phi
_j}^{-1}\Phi_j(X_j\,{-}\,\tilde{X}_J\setminus X_j)\bigr)\!\Biggr).\vadjust{\goodbreak}
\end{eqnarray*}
After estimates of the $\sigma^2_j$ and any parameters in the $\Phi_j$
have been plugged in,
the results in (\ref{cond1}) and (\ref{cond2}) can be combined to
obtain a Gaussian estimated predictive distribution for $f(x)$
conditional on $f|_{X_J}$ with mean given by (\ref{multi-step2}).
For generating confidence intervals, the variance of the estimated
predictive distribution, conditional on $f|_{X_J}$, can be calculated
in a backwards recursive manner using (\ref{cond1}) and~(\ref{cond2}).
Once again, note that $A^{-1}b$ should be taken as shorthand for the
solution to the linear system $Ax=b$.


\section{Nominal and numeric error}\label{nominalnumeric}

Now, we develop some intuition for why the multi-step procedure can
improve accuracy in many situations in practice.
First, computed quantities, which are subject to floating point error,
are distinguished from the idealized quantities that could be obtained
if a~computer performed calculations with full accuracy.
Hereafter, computed quantities will be distinguished with a tilde, such
as $\tilde{y}$.
We introduce the following separation of error into \textit{nominal} and
\textit{numeric} portions:
%
\begin{eqnarray}\label{errorSeparation}
|f(x)-\tilde{\mathcal{P}}(x)|&=&|f(x)-\mathcal{P}(x)+\mathcal
{P}(x)-\tilde{\mathcal{P}}(x)|\nonumber\\[-8pt]\\[-8pt]
&\le&|f(x)-\mathcal{P}(x)|+|\mathcal{P}(x)-\tilde{\mathcal
{P}}(x)|.\nonumber
\end{eqnarray}
Note that the absolute values in inequality (\ref{errorSeparation}) can
be replaced with the norm of one's choosing.
It is necessary to account for both nominal and numeric error since the
trade-off between the two is very important.
In most situations, reducing one will increase the other.
The following proposition shows that the native space norm of the
nominal error is always reduced by the addition of new data sites.
Throughout, let $\mathcal{N}_{\Phi}(\Omega)$ denote the reproducing
kernel Hilbert space corresponding to the positive definite kernel
$\Phi
$, and let~\mbox{$\|\cdot\|_{\mathcal{N}_{\Phi}(\Omega)}$} denote the norm on
that space
\cite{aronszajn}.
%
\begin{prop}\label{nominalImprovementProp}
If $f\in\mathcal{N}_{\Phi}(\Omega)$ and $X_1\subseteq X_2$, then
\[
\|f-\mathcal{P}_2\|_{\mathcal{N}_{\Phi}(\Omega)}\le\|f-\mathcal
{P}_1\|
_{\mathcal{N}_{\Phi}(\Omega)},
\]
where $\mathcal{P}_1$ and $\mathcal{P}_2$ denote the single-stage
interpolators on the sets $X_1$ and~$X_2$, respectively.
\end{prop}
\begin{pf}
It can be shown that the interpolator is orthogonal to its error with
respect to the native space inner product.
This implies that the result holds if and only if
\begin{eqnarray*}
&&\|f\|^2_{\mathcal{N}_{\Phi}(\Omega)}-\|\mathcal{P}_2\|^2_{\mathcal
{N}_{\Phi}(\Omega)}\le\|f\|^2_{\mathcal{N}_{\Phi}(\Omega)}-\|
\mathcal
{P}_1\|^2_{\mathcal{N}_{\Phi}(\Omega)}\\
&&\quad\Longleftrightarrow\quad\|\mathcal{P}_2\|^2_{\mathcal{N}_{\Phi}(\Omega
)}\ge
\|\mathcal{P}_1\|^2_{\mathcal{N}_{\Phi}(\Omega)}\\
&&\quad\Longleftrightarrow\quad f|'_{X_2}A^{-1}_{X_2,\Phi}f|_{X_2}\ge
f|'_{X_1}A^{-1}_{X_1,\Phi}f|_{X_1},
\end{eqnarray*}
where the last equivalent condition follows from the definition of the
native space norm\vadjust{\goodbreak} and the fact that $\alpha^j=A^{-1}_{X_j,\Phi
}f|_{X_j}$ for $A_{X_j,\Phi}=\{\Phi(x_u-x_v)\}$, $x_u,x_v\in X_j$, $j=1,2$.
Then, write the interpolation matrix $A_{X_2,\Phi}$ as
\[
A_{X_2,\Phi}=\pmatrix{
A_{X_1,\Phi} & A_{12}\cr
A_{21} & A_{22}},
\]
where $A_{12}=\Phi(X_1-X_2\setminus X_1)$, $A_{21}=\Phi(X_2\setminus
X_1-X_1)$, and $A_{22}=\Phi(X_2\setminus X_1-X_2\setminus X_1)$,
using the notation in (\ref{notation}).
Using partitioned matrix inverse and binomial inverse results
\cite{harville},
it can be shown that
\begin{eqnarray*}
&&f|'_{X_2}A^{-1}_{X_2,\Phi}f|_{X_2}\\
&&\qquad=f|'_{X_1}A^{-1}_{X_1,\Phi
}f|_{X_1}\\
&&\qquad\quad{} +
(f|_{X_2\setminus X_1}-A_{21}A^{-1}_{X_1,\Phi
}f|_{X_1})'A^{-1}_{22\cdot
1}(f|_{X_2\setminus X_1}-A_{21}A^{-1}_{X_1,\Phi}f|_{X_1}),
\end{eqnarray*}
where $A_{22\cdot1}=A_{22}-A_{21}A^{-1}_{X_1,\Phi}A_{12}$. Since
$A^{-1}_{22\cdot1}$ is a block on the diagonal of~$A^{-1}_{X_2,\Phi}$,
it must be positive definite and the result follows.
\end{pf}

On the other hand, the numeric error can become arbitrarily large by
the addition of new data sites.
Throughout, let $\lambda_{\max}(A)$ and $\lambda_{\min}(A)$
denote the maximum and minimum eigenvalues, respectively, of a positive
definite matrix $A$.
Note that $\lambda_{\min}(A_{X,\Phi})\to0$ as ${\min_{x_u\ne
x_v}}\|
x_u-x_v\|_2\to0$. Therefore, $\lambda_{\max}(A^{-1}_{X,\Phi})\to
\infty$ as ${\min_{x_u\ne x_v}}\|x_u-x_v\|_2\to0$.
An unboundedly large maximum eigenvalue of $A^{-1}_{X,\Phi}$ can
enormously amplify small errors in the function and kernel evaluations.
Consider the numeric error of the interpolator at a~new point $x$,
\begin{eqnarray*}
\mathcal{P}(x)-\tilde{\mathcal{P}}(x)&=&\sum_{i=1}^n[\alpha_i\Phi
(x-x_i)-\tilde{\alpha}_i\tilde{\Phi}(x-x_i)]\\[-2pt]
&=&\sum_{i=1}^n\bigl[(\alpha_i-\tilde{\alpha}_i)\Phi(x-x_i)-\tilde
{\alpha
}_i\bigl(\tilde{\Phi}(x-x_i)-\Phi(x-x_i)\bigr)\bigr].
\end{eqnarray*}
Let $\varepsilon^\alpha=\alpha-\tilde{\alpha}$ and $\varepsilon
^\Phi
=\tilde{\Phi}(X-x)-\Phi(X-x)$ using the notation in (\ref{notation}).
Then
\begin{eqnarray*}
\mathcal{P}(x)-\tilde{\mathcal{P}}(x)&=&\sum_{i=1}^n[\varepsilon
^\alpha
_i\Phi(x-x_i)-(\alpha_i-\varepsilon^\alpha_i)\varepsilon^\Phi_i]\\[-2pt]
&=&\sum_{i=1}^n[\varepsilon^\alpha_i\Phi(x-x_i)-\alpha_i\varepsilon
^\Phi
_i+\varepsilon^\alpha_i\varepsilon^\Phi_i].
\end{eqnarray*}
So,
%
\begin{equation}\label{lowerBound}
|\mathcal{P}(x)-\tilde{\mathcal{P}}(x)|
\ge|f|'_XA^{-1}_{X,\Phi}\varepsilon^\Phi|-\|\varepsilon^\alpha\|
_2\|\Phi
(x-X)\|_2-\|\varepsilon^\alpha\|_2\|\varepsilon^\Phi\|_2
\end{equation}
since $A_{X,\Phi}\alpha=f|_X$.
If, for example, $\varepsilon^\Phi$ is proportional to the eigenvector
corresponding to $\lambda_{\min}(A_{X,\Phi})$, and $f|_X$ is not
orthogonal to $\varepsilon^\Phi$, then the right-hand side of (\ref
{lowerBound}) can be made unboundedly large by taking
$\lambda_{\min}(A_{X,\Phi})\to0$.\vadjust{\goodbreak}

\begin{figure}

\includegraphics{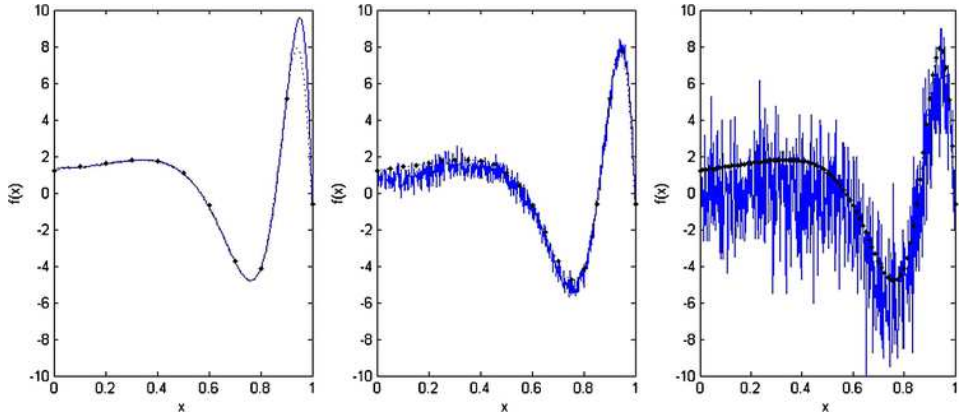}

\caption{Panels 1--3: interpolator in solid blue and actual
function in dotted black with collected data indicated by black dots.}
\label{numErrorExample}
\end{figure}

This phenomenon can be illustrated by attempting to build an
interpolator for the function
\[
f(x)=\operatorname{exp}\{(x+1/2)^2\}\operatorname{sin}\bigl(\operatorname{exp}
\{
(x+1/2)^2\}\bigr)
\]
shown in Figure \ref{numErrorExample} using the Gaussian kernel
\[
\Phi(x-y)=\operatorname{exp}\{-(x-y)^2\}.
\]
Interpolators, shown in blue, are built on 11, 21 and 81 evenly spaced
data points, shown in black dots, in the respective panels of Figure
\ref{numErrorExample}.
As the density of points increases, so does the numeric error.

%


Suppose that one is in the situation where most of the data sites are
well spread, but a few poorly separated data sites are causing small
numeric errors to be amplified. Consider forming an interpolator in two
stages. In the first stage, remove the data sites which are causing the
ill-conditioning of the interpolation matrix and interpolate the
remaining points with a~relatively wide kernel. The nominal error will
be only slightly larger than the error for the full data set, since the
removed data sites were nearly equal to data sites which were included.
However, the numeric error will be substantially less than that of an
interpolator formed on the full data set. In the second stage,
interpolate the residuals from the first-stage interpolator using a
kernel which is narrow enough that numeric errors remain small. The
second-stage interpolator will increase neither the nominal accuracy
nor the numeric error substantially. When the two interpolators are
added together to form the multi-step interpolator, the nominal
accuracy may be slightly worse, but the numeric accuracy will be very
much better.



\begin{figure}

\includegraphics{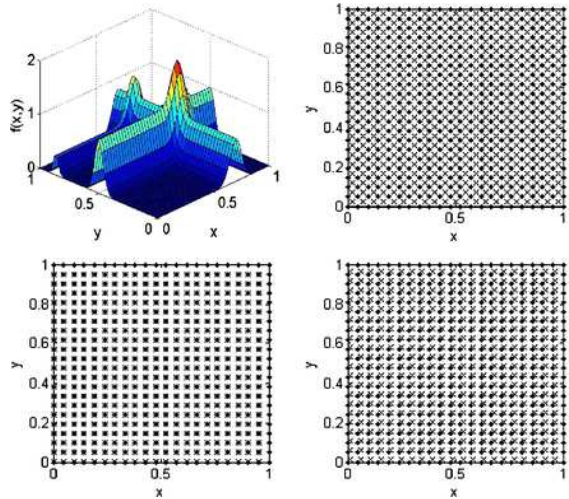}

\caption{Panel 1: the Michalewicz function. Panels 2--4 in
clockwise order: 925 point data sets with separation distances
$0.017$, $0.009$ and $5\times10^{-11}$, respectively.}\label{toyexample}
\end{figure}

For example, consider building an emulator for the Michalewicz function
\[
f(x,y)=\sin(\pi x)\sin^{20}(\pi x^2)+\sin(\pi y)
\sin^{20}(2\pi y^2)
\]
using the third 925 point data set in Figure \ref{toyexample} with
separation distance $5\times10^{-11}$.
%
The separation distance of a point set $X$ is half the distance between
the closest two points,
%
\begin{equation}\label{eqndistance}
q_X={\frac{1}{2}\min_{x_i,x_j\in X}}\|x_i-x_j\|_2.
\end{equation}
Clearly, the $\times$'s do not contribute much information about the
unknown surface.
If an ordinary Gaussian kernel interpolator, corresponding to a single
stage with
%
\begin{equation}\label{gaussianKernel}
\Phi(x-y)=\operatorname{exp}\Biggl\{-\sum_{j=1}^2\theta_j(x_j-y_j)^2\Biggr\}
\end{equation}
is built using all the data sites, the best possible mean squared
prediction error over values of $\theta_1,\theta_2$ is \mbox{$\approx$}$0.15$,
the square of the function's $L_2$ norm.
This is because the kernel must be very narrow, or the interpolation
matrix will be nearly singular.
Throughout, the term \textit{mean squared prediction error} is taken to
be the average prediction error over the input domain.
If, on the other hand, the $\bolds{\cdot}$'s are interpolated and
then the residuals on the $\times$'s are interpolated, corresponding to
two stages, the best possible mean squared prediction error over values
of $\theta_1,\theta_2$ at each stage is \mbox{$\approx$}$1.5\times10^{-5}$.
%
%
%
%

\section{Numeric accuracy}\label{secnum}

The numeric accuracy of the multi-step interpola\-tion procedure depends
on the accuracy of floating point matrix manipula\-tions.
Floating point accuracy refers to the fact that computers do not~%
per\-form calculations with real numbers, but instead with rounded
versions thereof.
For example, a typical computer has 15 digits of accuracy
meaning that
\[
\frac{\|\tilde{x}-x\|_2}{\|x\|_2}\le10^{-15},
\]
where $x$ denotes the actual value, and $\tilde{x}$ denotes the value
that the computer stores.

\subsection{Numeric accuracy of matrix inversion}\label{numSec1}

The following lemma on the accuracy of floating point matrix inversion
is a combination and generalization of results in
\cite{golub}.
%
\begin{definition}
$\!\!\!$The matrix 2-norm \mbox{$\|\cdot\|_2$} is defined as $\|A\|_2\,{=}\,\sqrt{\lambda
_{\max}(A'A)}$.
\end{definition}

\begin{lemma}\label{floatinglemma}
Suppose $Ax=b$ and $\tilde{A}\tilde{x}=\tilde{b}$ with $\|A-\tilde
{A}\|
_2\le\delta_A\|A\|_2$, $\|b-\tilde{b}\|_2\le\delta_b\|b\|_2$ and
$\kappa
(A)=r/\delta_A<1/\delta_A$ for some $\delta_A,\delta_b>0$. Then
$\tilde
{A}$ is nonsingular,
%
\begin{eqnarray}\label{floatingresult}
\frac{\|\tilde{x}\|_2}{\|x\|_2}&\le&\frac{1+r(\delta_b/\delta
_A)}{1-r},\nonumber\\[-9pt]\\[-9pt]
\frac{\|x-\tilde{x}\|_2}{\|x\|_2}&\le&\frac{\delta_A+\delta
_b}{1-r}\kappa
(A),\nonumber
\end{eqnarray}
where $\kappa(A)=\|A\|_2\|A^{-1}\|_2$.
\end{lemma}
\begin{pf}
Suppose $\tilde{A}$ is singular. Then there is a $y\ne0$ with $\tilde
{A}y=0$ so $(I-A^{-1}\tilde{A})y=y$. This implies $\|I-A^{-1}\tilde
{A}\|
_2\ge1$. On the other hand, the conditions $\|A-\tilde{A}\|_2\le
\delta
_A\|A\|_2$ and $\kappa(A)<1/\delta_A$ imply
$\|I-A^{-1}\tilde{A}\|_2<1$ giving a contradiction.

Now,\vspace*{1pt} $\tilde{A}\tilde{x}=\tilde{b}$ implies $A^{-1}\tilde{A}\tilde
{x}=A^{-1}(b-(b-\tilde{b}))=x+A^{-1}(\tilde{b}-b)$. The condition $\|
I-A^{-1}\tilde{A}\|_2\le r$ implies $\|A^{-1}\tilde{A}\|_2\ge1-r$ and
in turn
\begin{eqnarray*}
\|\tilde{x}\|_2&\le&\frac{1}{1-r}(\|x\|_2+\|A^{-1}\|_2\|\tilde{b}-b\|
_2)\\[-2pt]
&\le&\frac{1}{1-r}(\|x\|_2+\delta_b\|A^{-1}\|_2\|b\|_2)\\[-2pt]
&\le&\frac{1}{1-r}\biggl(\|x\|_2+r\frac{\delta_b\|b\|_2}{\delta_A\|A\|
_2}\biggr)\\[-2pt]
&\le&\frac{1}{1-r}\bigl(\|x\|_2+r(\delta_b/\delta_A)\|x\|_2\bigr),
\end{eqnarray*}
where the first inequality follows from the stated condition, the
triangle inequality, and the fact that $\|By\|_2\le\|B\|_2\|y\|_2$, the
second inequality follows from the condition $\|b-\tilde{b}\|_2\le
\delta
_b\|b\|_2$, the third inequality follows from the condition $\kappa
(A)=r/\delta_A$ and the final inequality follows from $\|b\|_2\le\|A\|
_2\|x\|_2$. Dividing by $\|x\|_2$ gives the first inequality in (\ref
{floatingresult}).\vadjust{\goodbreak}

Note that $A(\tilde{x}-x)=\tilde{b}-b-(\tilde{A}-A)\tilde{x}$. So,
\begin{eqnarray*}
\|\tilde{x}-x\|_2&\le&\|A^{-1}\|_2\|\tilde{b}-b\|_2+\|A^{-1}\|_2\|
\tilde
{A}-A\|_2\|\tilde{x}\|_2\\
&\le&\delta_b\|A^{-1}\|_2\|b\|_2+\delta_A\|A^{-1}\|_2\|A\|_2\|\tilde
{x}\|
_2\\
&\le&\delta_b\kappa(A)\frac{\|b\|_2}{\|A\|_2}+\delta_A\kappa(A)\|
\tilde
{x}\|_2\\
&\le&\kappa(A)\|x\|_2\biggl(\delta_b+\delta_A\frac{1+r(\delta
_b/\delta
_A)}{1-r}\biggr),
\end{eqnarray*}
where the first inequality follows from the triangle inequality and the
fact that $\|By\|_2\le\|B\|_2\|y\|_2$, the second inequality follows
from the conditions $\|b-\tilde{b}\|_2\le\delta_b\|b\|_2$ and $\|
A-\tilde{A}\|_2\le\delta_A\|A\|_2$, the third inequality follows from
the definition of $\kappa(A)$ and the final inequality follows from the
fact that $\|b\|_2\le\|A\|_2\|x\|_2$ and the first inequality in (\ref
{floatingresult}). Dividing by $\|x\|_2$ and simplifying gives the
second part of (\ref{floatingresult}).
\end{pf}

\subsection{Numeric accuracy of single-stage interpolator}\label{numSec2}

The above lemma can be used to bound the numeric error of an
interpolator as follows.
%
\begin{theorem}\label{numericResult}
Suppose that $\|A_{X,\Phi}-\tilde{A}_{X,\Phi}\|_2\le\delta_A\|
A_{X,\Phi
}\|_2$, $\|f|_X-\break\tilde{f}|_X\|_2\le\delta_f\|f|_X\|_2$, $\kappa
(A_{X,\Phi})=r/\delta_A<1/\delta_A$ and $\sup_{x,y\in\Omega
}|\Phi
(x-y)-\tilde{\Phi}(x-y)|<D\delta_A$ for some $\delta_A,\delta
_f,D>0$, then
\begin{eqnarray*}
|\mathcal{P}(x)-\tilde{\mathcal{P}}(x)|&\le&
\bigl\|f|_X/\sqrt{n}\bigr\|_2\frac{(\delta_A+\delta_f)}{1-r}g(X,\Phi),\\
g(X,\Phi)&=&\frac{n}{\lambda_{\min}(A_{X,\Phi})}\bigl(\kappa
(A_{X,\Phi
})\Phi(0)+D\bigr),
\end{eqnarray*}
where $\kappa(\cdot)$ is defined in Lemma \ref{floatinglemma}.
\end{theorem}

Note that for large $n$ and approximately uniform $X$, $\|
f|_X/\sqrt{n}\|_2\approx\|f\|_{L_2(\Omega)}$, where
\[
\|f\|_{L_2(\Omega)}=\sqrt{\int_\Omega f(x)^2\,\mathrm{d}x}.
\]
Further, the assumption ${\sup_{x,y\in\Omega}}|\Phi(x-y)-\tilde{\Phi
}(x-y)|<D\delta_A$ requires that the kernel is computed in a relatively
accurate manner.
\begin{pf*}{Proof of Theorem \ref{numericResult}}
First,
\begin{eqnarray*}
\mathcal{P}(x)-\tilde{\mathcal{P}}(x)&=&\sum_{i=1}^n[\alpha_i\Phi
(x-x_i)-\tilde{\alpha}_i\tilde{\Phi}(x-x_i)]\\
&=&\sum_{i=1}^n\bigl[(\alpha_i-\tilde{\alpha}_i)\Phi(x-x_i)-\tilde
{\alpha
}_i\bigl(\tilde{\Phi}(x-x_i)-\Phi(x-x_i)\bigr)\bigr].
\end{eqnarray*}
So,
\begin{eqnarray*}
&&|\mathcal{P}(x)-\tilde{\mathcal{P}}(x)|\\
&&\qquad\le\Biggl|\sum_{i=1}^n(\alpha_i-\tilde{\alpha}_i)\Phi(x-x_i)\Biggr|\\
&&\qquad\quad{}
+\Biggl|\sum
_{i=1}^n\tilde{\alpha}_i\bigl(\tilde{\Phi}(x-x_i)-\Phi(x-x_i)\bigr)\Biggr|.
\end{eqnarray*}
Applying the Cauchy--Schwarz inequality to each term gives
\begin{eqnarray*}
%
&&|\mathcal{P}(x)-\tilde{\mathcal{P}}(x)|\\
&&\qquad\le\|\alpha-\tilde{\alpha}\|_2\sqrt{\sum_{i=1}^n\Phi(x-x_i)^2}\\
&&\qquad\quad{} +\|\tilde{\alpha}\|_2\sqrt{\sum_{i=1}^n\bigl(\tilde{\Phi
}(x-x_i)-\Phi
(x-x_i)\bigr)^2}.
\end{eqnarray*}
The terms under the radicals can be bounded to obtain
\[
|\mathcal{P}(x)-\tilde{\mathcal{P}}(x)| \le\sqrt{n}\|\alpha
-\tilde
{\alpha}\|_2\Phi(0)+\sqrt{n}\|\tilde{\alpha}\|_2 D\delta_A.
\]
Now, Lemma \ref{floatinglemma} can be applied to the coefficients, giving
\begin{eqnarray*}
|\mathcal{P}(x)-\tilde{\mathcal{P}}(x)| &\le& \sqrt{n}\frac{\delta
_A+\delta_f}{1-r}\kappa(A_{X,\Phi})\|\alpha\|_2\Phi(0)\\
&&{} +\sqrt{n}\frac{1+r(\delta_f/\delta_A)}{1-r}\|\alpha\|_2
D\delta
_A.
%
\end{eqnarray*}
Noting that $\|\alpha\|_2\le\|A^{-1}_{X,\Phi}\|_2\|f|_X\|_2$ and
collecting terms shows that
\begin{eqnarray*}
|\mathcal{P}(x)-\tilde{\mathcal{P}}(x)| &\le& \frac{\sqrt{n}\|
A^{-1}_{X,\Phi}\|_2\|f|_X\|_2}{1-r}\\
&&{} \times\bigl((\delta
_A+\delta
_f)\kappa(A_{X,\Phi})\Phi(0)+D(\delta_A+r\delta_f)\bigr)\\
&\le&\frac{\sqrt{n}\|A^{-1}_{X,\Phi}\|_2\|f|_X\|_2}{1-r}(\delta
_A+\delta
_f)\bigl(\kappa(A_{X,\Phi})\Phi(0)+D\bigr).
\end{eqnarray*}
Rearranging gives the result.
\end{pf*}

\subsection{Numeric accuracy of multi-step interpolator}\label{numSec3}

The first numeric result for the multi-step interpolator follows from
Theorem \ref{numericResult}.
Here, $\delta$ denotes the computer's floating point accuracy,
typically $\delta\le10^{-15}$.
%
\begin{theorem}\label{numericThm}
Suppose that for $j=1,\ldots,J$, $\|A_{X_j,\Phi_j}-\tilde
{A}_{X_j,\Phi
_j}\|_2\le\break\delta_j\|A_{X_j,\Phi_j}\|_2$, $\|f|_{X_j}-\tilde
{f}|_{X_j}\|
_2\le\delta\|f|_{X_j}\|_2$, $\kappa(A_{X_j,\Phi_j})\le r/\delta
_j<1/\delta_j$ and\break ${\sup_{x,y\in\Omega}}|\Phi_{j}(x-y)-\tilde
{\Phi
}_{j}(x-y)|<D\delta$ for some $\delta_j,\delta,D>0$ with $\delta_j\|
(f-\sum_{k=1}^{j-1}\mathcal{P}^k)|_{X_j}/\sqrt{n_j}\|_2\le\delta\|
f|_{X_j}/\sqrt{n_j}\|_2$, then
%
\begin{eqnarray}\label{numericThmEqn}\qquad
&&\Biggl|\sum_{j=1}^J\mathcal{P}^j(x)-\sum_{j=1}^J\tilde{\mathcal
{P}}^j(x)\Biggr|\nonumber\\[-9pt]\\[-9pt]
&&\qquad\le\delta\bigl\|f|_{X_J}/\sqrt{n_J}\bigr\|_2\Biggl[\sum_{M=1}^J C^M\sum
_{i\in
\mathcal{S}_J(M)} \prod_{k=1}^M\rho
(X_{i_k},X_{i_{k+1}})g(X_{i_k},\Phi
_{i_k})\Biggr],\nonumber
\end{eqnarray}
where $C=2/(1-r)$, $\mathcal{S}_J(M)=\{i\in\mathbb{N}^{M+1}\dvtx1\le
i_1<\cdots<i_M\le i_{M+1}=J\}$ $\rho(X,Y)=\|f|_X/\sqrt{n_X}\|_2/\|
f|_Y/\sqrt{n_Y}\|_2$, and $g$ is defined in Theorem \ref{numericResult}.
\end{theorem}

The assumption $\delta_j\|(f-\sum_{k=1}^{j-1}\mathcal
{P}^k)|_{X_j}/\sqrt{n_j}\|_2\le\delta\|f|_{X_j}/\sqrt{n_j}\|_2$ roughly
requires that the nominal errors either shrink or are not much larger
than the function values.
In practice, combinations of functions and training data sets which do
not meet this assumption are {very} rare.
\begin{pf*}{Proof of Theorem \ref{numericThm}}
The result can be shown using induction on the number of stages $J$. If
$J=1$, then the result follows immediately from Theorem~\ref{numericResult}.
Take $J\ge2$, and assume the result holds for $J-1$ stages. Then
%
\begin{eqnarray}\label{inequalities}
&&\Biggl\|\Biggl(f-\sum_{j=1}^{J-1}\mathcal{P}^j\Biggr)\Bigg|_{X_J}-\Biggl(\tilde{f}-\sum
_{j=1}^{J-1}\tilde{\mathcal{P}}^j\Biggr)\Bigg|_{X_J}\Biggr\|_2\nonumber\\[-2pt]
&&\qquad\le\bigl\|f|_{X_J}-\tilde{f}|_{X_J}\bigr\|_2+\Biggl\|\Biggl(\sum_{j=1}^{J-1}\mathcal
{P}^j-\sum_{j=1}^{J-1}\tilde{\mathcal{P}}^j\Biggr)\Bigg|_{X_J}\Biggr\|_2\\[-2pt]
&&\qquad\le\delta\|f|_{X_J}\|_2+\sqrt{n_J}\Biggl\|\sum_{j=1}^{J-1}\mathcal
{P}^j-\sum
_{j=1}^{J-1}\tilde{\mathcal{P}}^j\Biggr\|_{L_\infty(\Omega)},\nonumber
\end{eqnarray}
where the first inequality follows from the triangle inequality, and
the second inequality follows from the assumptions and by bounding the
$L_2$ error with the maximum error. The induction hypothesis can be
applied to the final term in (\ref{inequalities}) giving the bound
%
\begin{eqnarray}\label{bound}
&&\delta\|f|_{X_J}\|_2\nonumber\hspace*{-25pt}\\[-9pt]\\[-9pt]
&&\quad{}\times\Biggl(1\,{+}\,\rho(X_{J-1},X_J)\!\sum_{M=1}^{J-1}
C^M\!\sum
_{i\in\mathcal{S}_{J-1}(M)}\!\prod_{k=1}^M\!\rho
(X_{i_k},X_{i_{k+1}})g(X_{i_k},\Phi_{i_k})\!\Biggr).\hspace*{-25pt}\nonumber
\end{eqnarray}
In stage\vspace*{1pt} $J$, the error from the first $J-1$ stages are interpolated on
$X_J$. After multiplying\vadjust{\goodbreak} and dividing the above bound (\ref{bound}) by
$\|(f-\sum_{j=1}^{J-1}\mathcal{P}^j)|_{X_J}\|_2$, Theorem~\ref
{numericResult} can be used to bound the error due to stage $J$.
Note that the term $\delta_f$ in Theorem~\ref{numericResult} is the
above\vspace*{2pt} bound (\ref{bound}) divided by $\|(f-\sum_{j=1}^{J-1}\mathcal
{P}^j)|_{X_J}\|_2$ and the term $\delta_A$ in Theorem \ref
{numericResult} is $\delta_j$.
By assumption, $\delta_j$ is smaller than or equal to (\ref{bound})
divided by $\|(f-\sum_{j=1}^{J-1}\mathcal{P}^j)|_{X_J}\|_2$.
Simplification and coarsening of the bound gives
%
\begin{eqnarray}\label{bound2}
&&|\mathcal{P}^J(x)-\tilde{\mathcal{P}}^J(x)|\nonumber\hspace*{-25pt}\\
&&\quad\le\delta\bigl\|f|_{X_J}/\sqrt{n_J}\bigr\|_2\frac{2}{1-r}g(X_J,\Phi_J)\hspace*{-25pt}\\
&&\qquad\times\Biggl(1\,{+}\,\rho(X_{J-1},X_J)\!\sum_{M=1}^{J-1} C^M\!\sum
_{i\in\mathcal{S}_{J-1}(M)}\!\prod_{k=1}^M\rho
(X_{i_k},X_{i_{k+1}})g(X_{i_k},\Phi_{i_k})\!\Biggr).\hspace*{-25pt}\nonumber
\end{eqnarray}
Now,
\[
\Biggl|\sum_{j=1}^J\mathcal{P}^j(x)-\sum_{j=1}^J\tilde{\mathcal
{P}}^j(x)\Biggr| \le
\Biggl|\sum_{j=1}^{J-1}\mathcal{P}^j(x)-\sum_{j=1}^{J-1}\tilde{\mathcal
{P}}^j(x)\Biggr|+|\mathcal{P}^J(x)-\tilde{\mathcal{P}}^J(x)|.
\]
So, the induction hypothesis can be applied again along with (\ref
{bound2}) giving
%
\begin{eqnarray}\label{bound3}
&&\Biggl|\sum_{j=1}^J\mathcal{P}^j(x)-\sum_{j=1}^J\tilde{\mathcal
{P}}^j(x)\Biggr|\nonumber\\
&&\qquad\le\delta\bigl\|f|_{X_J}/\sqrt{n_J}\bigr\|_2\nonumber\\
&&\qquad\quad{} \times\Biggl[\rho(X_{J-1},X_J)\sum_{M=1}^{J-1} C^M\sum_{i\in
\mathcal{S}_{J-1}(M)} \prod_{k=1}^M\rho
(X_{i_k},X_{i_{k+1}})g(X_{i_k},\Phi_{i_k})\nonumber\\[-8pt]\\[-8pt]
&&\qquad\quad\hspace*{17.5pt}{} + C g(X_J,\Phi_J)\nonumber\\
&&\qquad\quad\hspace*{17.5pt}{} + C \rho(X_{J-1},X_J)g(X_J,\Phi_J) \nonumber\\
&&\qquad\quad\hspace*{68.2pt}{}\times\sum
_{M=1}^{J-1} C^M\sum_{i\in\mathcal{S}_{J-1}(M)} \prod_{k=1}^M\rho
(X_{i_k},X_{i_{k+1}})g(X_{i_k},\Phi_{i_k})\Biggr].\nonumber
\end{eqnarray}
Note that the term in square brackets in (\ref{numericThmEqn}) is the
sum of the terms with $i_M<J$ and $i_M=J$ giving
\begin{eqnarray*}
&&\sum_{M=1}^J C^M\sum_{i\in\mathcal{S}_J(M)} \prod_{k=1}^M\rho
(X_{i_k},X_{i_{k+1}})g(X_{i_k},\Phi_{i_k})\\
&&\qquad=\sum_{M=1}^{J-1} C^M\sum_{i\in\mathcal{S}_J(M),i_M<J} \prod
_{k=1}^M\rho(X_{i_k},X_{i_{k+1}})g(X_{i_k},\Phi_{i_k})\\
&&\qquad\quad{} +\sum_{M=1}^J C^M\sum_{i\in\mathcal{S}_J(M),i_M=J} \prod
_{k=1}^M\rho(X_{i_k},X_{i_{k+1}})g(X_{i_k},\Phi_{i_k})\\
&&\qquad=\rho(X_{J-1},X_J)\sum_{M=1}^{J-1} C^M\sum_{i\in\mathcal
{S}_{J-1}(M)}
\prod_{k=1}^M\rho(X_{i_k},X_{i_{k+1}})g(X_{i_k},\Phi_{i_k})\\
&&\qquad\quad{}+Cg(X_j,\Phi_J)+\sum_{M=2}^J C^M\sum_{i\in\mathcal
{S}_J(M),i_M=J} \prod_{k=1}^M\rho(X_{i_k},X_{i_{k+1}})g(X_{i_k},\Phi
_{i_k}),
\end{eqnarray*}
which is exactly the term in square brackets in (\ref{bound3}), proving
the result.
\end{pf*}

\subsection{Dependence on separation distance}

The terms
%
\begin{equation}\label{g}
g(X_j,\Phi_j)=\frac{n_j}{\lambda_{\min}(A_{X_j,\Phi_j})}\bigl(\kappa
(A_{X_j,\Phi_j})\Phi(0)+D\bigr)
\end{equation}
from Theorem \ref{numericThm} can be computed, at least approximately.
However, by bounding (\ref{g}) in terms of the separation distance, as
defined in (\ref{eqndistance}),
the role of the data sites and the kernel's smoothness in the numeric
accuracy are revealed.
These results indicate that using poorly separated data or a wide
kernel $\Phi$ with a rapidly decaying Fourier transform, implying more
smoothness, has more potential to result in large numeric errors in
interpolation.
The Fourier transform can be defined as follows.
%
\begin{definition}
For $f\in L_1(\mathbb{R}^d)$ define the Fourier transform \cite{stein1971}
\[
\hat{f}(\omega)=(2\pi)^{-d/2}\int_{\mathbb{R}^d}f(x)e^{-i\omega
'x}\,\mathrm{d}x.
\]
\end{definition}

%
To generate the bound on (\ref{g}), the following result from
\cite{wendland2005}
can be used.
%
\begin{theorem}\label{mineigenthm}
Let $\varphi_*(M,\Phi)=\inf_{\|\omega\|_2\le2M}\hat{\Phi}(\omega
)$. Then
\begin{eqnarray*}
\lambda_{\min}(A_{X,\Phi})&\ge& C_d\varphi_*(M_d/q,\Phi)/q^d,\\
M_d&=&12\bigl(\pi\Gamma^2(d/2+1)/9\bigr)^{1/(d+1)},\\
C_d&=&(M_d/2^{3/2})^d/\bigl(2\Gamma(d/2+1)\bigr)
\end{eqnarray*}
for any $q\le q_X$, where $A_{X,\Phi}=\{\Phi(x_i-x_j)\}$.
\end{theorem}

To bound $\lambda_{\max}(A_{X,\Phi})$ below, Gershgorin's theorem
\cite{varga} can be used.
Gershgorin's theorem
states that the largest eigenvalue of $A_{X,\Phi}$ has
\[
|\lambda_{\max}(A_{X,\Phi})-\Phi(x_j-x_j)|\le\sum_{i=1,i\ne
j}^{n}|\Phi(x_i-x_j)|.
\]
Rearranging and coarsening the bound gives
%
\begin{equation}\label{maxeigenbound}
\lambda_{\max}(A_{X,\Phi})\le n\Phi(0).
\end{equation}
Theorem \ref{mineigenthm} and inequality (\ref{maxeigenbound}) can be
combined to obtain the following theorem bounding (\ref{g}).
%
\begin{theorem}\label{thmg}
Under the assumptions in Theorem \ref{numericThm},
\begin{eqnarray*}
g(X_j,\Phi_j)&\le&\kappa_{\mathrm{upper}}(X_j,\Phi_j)\bigl(\kappa_{
\mathrm{upper}}(X_j,\Phi_j)\Phi(0)+D\bigr),\\
\kappa_{\mathrm{upper}}(X_j,\Phi_j)&=&\frac{n_j q^d_{X_j}}{C_d\varphi
_*(M_d/q_{X_j},\Phi_j)}.
\end{eqnarray*}
\end{theorem}

The nested sequence $X_1\subset\cdots\subset X_J$ in (\ref{eqnnested-design})
with large separation distance can be generated from \textit{nested
space-filling designs} \cite
{qian20091,qtw2009,qa10,qian20092,haaland1}, which
were originally developed for the purpose of conducting multi-fidelity
computer experiments.
Space-filling designs have shown particular merit in numerical
integration \cite
{stein1987,Owen1992,owen1994,owen1997,OWEN1997a,tang93,loh96a,loh96b,loh08,loh}.
Theorem \ref{thmg} provides new
insights into the use of such designs
in interpolation.

\section{Nominal accuracy}\label{nominalsection}

The results in this section indicate that the nominal error in
interpolation converges to zero more quickly for wider, smoother
kernels~$\Phi$, although the constant involved in this rate changes.
This is in direct opposition to the numeric error, which tends to be
smaller for narrower, less smooth kernels.
In fact, it will be seen that convergence of the nominal error of an
arbitrarily fast \textit{rate} can be achieved with an infinitely smooth
kernel, such as the Gaussian in (\ref{gaussianKernel}).

A re-scaling is introduced in the following definition.
%
\begin{definition}
For a nonsingular $\Theta$, define $\Phi_\Theta(x)=\Phi(\Theta x)$.
\end{definition}

\subsection{Point-wise bound}

Initially, consider a single stage with a fixed $\Phi$ which is
re-scaled by a fixed $\Theta$. For a set of input sites $X$ of size
$n$, define the cardinal basis functions
\begin{eqnarray*}
u_i(x)&=&\sum_{i=1}^n\beta_i\Phi_\Theta(x-x_j),\\
u_i(x_j)&=&\ind_{\{i=j\}}
\end{eqnarray*}
for $i,j=1,\ldots,n$. Then
\[
\mathcal{P}(x)=\sum_{i=1}^n f(x_i)u_i(x).
\]
Since $f(x)=\langle f,\Phi_\Theta(\cdot-x)\rangle_{\mathcal
{N}_{\Phi
_\Theta}(\Omega)}$
if $f\in\mathcal{N}_{\Phi_\Theta}(\Omega)$,
\begin{eqnarray*}
f(x)-\mathcal{P}(x)&=&\langle f,\Phi_\Theta(\cdot-x)\rangle
_{\mathcal
{N}_{\Phi_\Theta}(\Omega)}-\sum_{i=1}^n u_i(x)\langle f,\Phi
_\Theta
(\cdot-x_i)\rangle_{\mathcal{N}_{\Phi_\Theta}(\Omega)}\\
&=&\Biggl\langle f,\Phi_\Theta(\cdot-x)-\sum_{i=1}^n u_i(x)\Phi_\Theta
(\cdot
-x_i)\Biggr\rangle_{\mathcal{N}_{\Phi_\Theta}(\Omega)}.%
\end{eqnarray*}
Now, the Cauchy--Schwarz inequality can be applied, giving the error bound
%
\begin{equation}\label{bound4}
|f(x)-\mathcal{P}(x)|\le\|f\|_{\mathcal{N}_{\Phi_\Theta}(\Omega
)}\Biggl\|\Phi
_\Theta(\cdot-x)-\sum_{i=1}^n u_i(x)\Phi_\Theta(\cdot-x_i)\Biggr\|
_{\mathcal
{N}_{\Phi_\Theta}(\Omega)}.\hspace*{-30pt}
\end{equation}
The second term on the right-hand side of (\ref{bound4}) is the
so-called \textit{power function}, $P_{\Phi_\Theta,X}$.
It can be shown
\cite{wendland2005}
that if the domain of interest $\Omega$ is bounded and convex, then
\[
P^2_{\Phi_\Theta,X}\le C_1\|\Phi_\Theta-p\|_{L_\infty(B(0,C_2 h_X))},
\]
where $C_1,C_2>0$ are constants which may depend on $\Omega$, $p$ is
any multivariate polynomial, $B(a,b)=\{x\in\mathbb{R}^d\dvtx\|x-a\|_2<b\}$
and $h_X$ denotes the \textit{fill distance}
\[
h_X={\sup_{x\in\Omega}\min_{x_u\in X}}\|x-x_u\|_2.
\]
Now, if $\Phi$ has $k$ continuous derivatives, $p$ can be taken to be
the Taylor's polynomial of $\Phi_\Theta$ of degree $k-1$. Then
\[
\|\Phi_\Theta-p\|_{L_\infty(B(0,C_2h_X))}\le C_3 \|\Theta\|^k_2 h^k_X,
\]
where $C_3$ is a constant which does not depend on $\Theta$. Combining
the above development gives the following.
%
\begin{theorem}\label{nominalResult}
Suppose that $\Omega$ is bounded and convex, $\Phi$ satisfies
Assumption \ref{kernelAssumption} and has $k$ continuous derivatives
and $\Theta$ is nonsingular. Then
\[
|f(x)-\mathcal{P}(x)|\le C_{\Phi} \|\Theta\|^{k/2}_2 h^{k/2}_X\|f\|
_{\mathcal{N}_{\Phi_\Theta}(\Omega)}.
\]
%
\end{theorem}

\subsection{Native space bound}

First, 
write $\Phi_{\Theta}\ast\Phi_{\Theta}$ as
\[
\Phi_{\Theta}\ast\Phi_{\Theta}(x-y)=\int_\Omega\Phi_{\Theta
}(x-t)\Phi
_{\Theta}(y-t)\,\mathrm{d}t.
\]
%
Then, for $f\in\mathcal{N}_{\Phi_{\Theta}\ast\Phi_{\Theta
}}(\Omega)$
and $x\in\Omega$, express $f$ in terms of the integral operator
\[
f(x)=\int_\Omega u(y)\Phi_{\Theta}\ast\Phi_{\Theta}(x-y)\,\mathrm{d}y,
\]
where $u\in L_2(\Omega)$. Combining these expressions gives
\begin{eqnarray*}
f(x)&=&\int_\Omega u(y)\int_\Omega\Phi_\Theta(y-t)\Phi_\Theta
(x-t)\,\mathrm{d}t\,\mathrm{d}y\\
&=&\int_\Omega v(t)\Phi_{\Theta}(x-t)\,\mathrm{d}t,
\end{eqnarray*}
where $v\in L_2(\Omega)$ is given by
\[
v(t)=\int_\Omega u(y)\Phi_{\Theta}(y-t)\,\mathrm{d}y
\]
for $t\in\Omega$. Then
%
\begin{eqnarray}\label{nominalbounddevelopment}
\|f-\mathcal{P}\|^2_{\mathcal{N}_{\Phi_{\Theta}}(\Omega)}&=&
\langle f-\mathcal{P},f\rangle_{\mathcal{N}_{\Phi_{\Theta}}(\Omega
)}\nonumber\\
&=&\langle f-\mathcal{P},v\rangle_{L_2(\Omega)}\\
&\le&\|f-\mathcal{P}\|_{L_2(\Omega)}\|v\|_{L_2(\Omega)},\nonumber
\end{eqnarray}
where the first equality follows from the orthogonality of the
interpolator and its error with respect to the native space norm, the
second equality follows from the properties of the integral operator
and the inequality follows from the Cauchy--Schwarz inequality.


If $\Phi$ has $k$ continuous derivatives, then the first term on the
right-hand side of inequality (\ref{nominalbounddevelopment}) can be
bounded using Theorem \ref{nominalResult} as
%
\begin{eqnarray}\label{L2bound}
\|f-\mathcal{P}\|_{L_2(\Omega)}&\le&\sqrt{\operatorname{vol} \Omega} \|
f-\mathcal{P}\|_{L_\infty}(\Omega)\\
&\le& C_{\Phi}\|\Theta\|^{k/2}_2h_{X}^{k/2}\|f-\mathcal{P}\|
_{\mathcal
{N}_{\Phi_\Theta}(\Omega)},\nonumber
\end{eqnarray}
where the first inequality follows by relating the $L_2(\Omega)$ and
$L_\infty(\Omega)$ norms, and the second inequality follows by applying
Theorem \ref{nominalResult} to $f-\mathcal{P}$. Plugging inequality
(\ref{L2bound}) into inequality (\ref{nominalbounddevelopment}) and
canceling a single $\|f-\mathcal{P}\|_{\mathcal{N}_{\Phi_{\Theta
}}(\Omega)}$ term gives
%
\begin{equation}\label{nominalbounddevelopment2}
\|f-\mathcal{P}\|_{\mathcal{N}_{\Phi_{\Theta}}(\Omega)}\le C_{\Phi
}\|
\Theta\|^{k/2}_2h_{X}^{k/2}\|v\|_{L_2(\Omega)}.
\end{equation}
Using the properties of the integral operator, the square of the second
term on the right-hand side of inequality (\ref
{nominalbounddevelopment2}) can be expressed as
%
\begin{eqnarray}\label{nominalbounddevelopment3}
\|v\|^2_{L_2(\Omega)}&=&\int_{\Omega^3}u(x)u(y)\Phi_{\Theta
}(y-t)\Phi
_{\Theta}(x-t)\,\mathrm{d}x\,\mathrm{d}y\,\mathrm{d}t\nonumber\\[-8pt]\\[-8pt]
&=&\|f\|^2_{\mathcal{N}_{\Phi_{\Theta}\ast\Phi_{\Theta}}(\Omega
)}.\nonumber
\end{eqnarray}
Combining inequality (\ref{nominalbounddevelopment2}) and equality
(\ref
{nominalbounddevelopment3}) gives the following theorem.

\begin{theorem}\label{nativeSpaceBound}
Under the assumptions of Theorem \ref{nominalResult},
\[
\|f-\mathcal{P}\|_{\mathcal{N}_{\Phi_\Theta}(\Omega)}\le C_\Phi\|
\Theta
\|_2^{k/2} h_{X}^{k/2}\|f\|_{\mathcal{N}_{\Phi_\Theta\ast\Phi
_\Theta
}(\Omega)}.\vadjust{\goodbreak}
\]
\end{theorem}

To allow for individual re-scalings in different stages, we start with
some notation. Define $\Psi_k$ recursively as
%
\begin{eqnarray}\label{kernelRelation1}
\Psi^0&=&\Phi,\nonumber\\[-8pt]\\[-8pt]
\Psi^k&=&\Psi^{k-1}\ast\Psi^{k-1}\nonumber
\end{eqnarray}
for $k\in\mathbb{N}$.
For the kernel on step $j$, take
%
\begin{equation}\label{kernelRelation2}
\Phi_j=\Psi^{J-j}_{\Theta_j}.
\end{equation}
We now develop a bound on \mbox{$\|\cdot\|_{\mathcal{N}_{\Phi_j\ast\Phi
_j}(\Omega)}$} in terms of \mbox{$\|\cdot\|_{\mathcal{N}_{\Phi
_{j-1}}(\Omega
)}$}. The \mbox{basic} assumptions on the re-scaling matrices $\Theta_j$ in
this section are that they are nonsingular and \textit{larger} than the
$\Theta_{j-1}$ in the sense that $\lambda_{\max}(\Theta
'_{j-1}\Theta
_{j-1}\Xi'_j\Xi_j)\le1$, where $\Xi'_j=\Theta_j^{-1}$.


%
%

In the case $\Omega=\mathbb{R}^d$, the native space $\mathcal
{N}_{\Phi
_\Theta}(\mathbb{R}^d)$
has norm defined through the inner product
%
\begin{equation}\label{InnerProdRd}
\langle f,g\rangle_{\mathcal{N}_{\Phi_\Theta}(\mathbb{R}^d)}=(2\pi
)^{-d/2} \int_{\mathbb{R}^d}\frac{\hat{f}(\omega)\overline{\hat
{g}(\omega)}}{\hat{\Phi}_\Theta(\omega)}\,\mathrm{d}\omega,
\end{equation}
where $\hat{f}$ and $\overline{\hat{g}}$ denote the Fourier transform
and complex conjugate of the Fourier transform, respectively, of
$f,g\in
\mathcal{N}_{\Phi_\Theta}(\mathbb{R}^d)$
\cite{wendland2005}.
This explicit representation of the native space inner product can be
used\vspace*{1pt} to relate the native space norms for convolutions and re-scalings.
Hereafter, take
$\infty>c_2\ge c_1>0$ and $\hat{\Upsilon}$ with
%
\begin{equation}\label{c1c2}
\omega'\omega\le\nu'\nu\quad\Longrightarrow\quad\hat{\Upsilon}(\omega
)\ge\hat
{\Upsilon}(\nu),\qquad
c_1\hat{\Upsilon}(\omega)\le\hat{\Phi}(\omega)\le c_2\hat
{\Upsilon
}(\omega).\hspace*{-15pt}
\end{equation}
Assumption \ref{kernelAssumption} ensures that $c_1$, $c_2$ and $\hat
{\Psi}$ satisfying (\ref{c1c2}) exist \cite{wendland2005}.
The bounds to follow are tightest for $c_2-c_1$ as small as possible.
Essentially, we want a~\textit{radially decreasing} envelop on the
Fourier transform of the underlying kernel~$\Phi$ to simplify development.
Note that the Fourier transforms $\hat{\Phi}$ and~$\hat{\Phi
}_\Theta$
are related in the following manner:
%
\begin{eqnarray}\label{FourierRelationship}
\hat{\Phi}_\Theta(\omega)&=&(2\pi)^{-d/2}\int_{\mathbb{R}^d}\Phi
_\Theta
(x)e^{-i\omega'x}\,\mathrm{d}x\nonumber\\
&=&(2\pi)^{-d/2}\int_{\mathbb{R}^d}\Phi(\Theta x)e^{-i\omega'\Xi
'\Theta
x}\,\mathrm{d}x\nonumber\\[-8pt]\\[-8pt]
&=&(2\pi)^{-d/2}|{\operatorname{det}}(\Xi)|\int_{\mathbb{R}^d}\Phi
(y)e^{-i\omega'\Xi
'y}\,\mathrm{d}y\nonumber\\
&=&|{\operatorname{det}}(\Xi)|\hat{\Phi}(\Xi\omega),\nonumber
\end{eqnarray}
where $\Xi'=\Theta^{-1}$ and the third equality follows by making the
substitution
$y=\Theta x$.
%
\begin{prop}\label{NormRelationProposition}
If Assumption \ref{kernelAssumption}
is satisfied and $\Theta_{j-1},\Theta_j$ are nonsingular with
respective inverses $\Xi'_{j-1},\Xi'_j$,
then
%
\begin{eqnarray}\label{NormRelationEquation}
&&\lambda_{\max}(\Theta'_{j-1}\Theta_{j-1}\Xi'_j\Xi_j)\le1\nonumber\\[-8pt]\\[-8pt]
&&\quad\Longrightarrow\quad\|f\|^2_{\mathcal{N}_{\Phi_j\ast\Phi_j}(\mathbb
{R}^d)}\le\biggl(\frac{c_2}{c_1}\biggr)^{2^{J-(j-1)}}\frac{|{\operatorname{det}}(\Xi_{j-1})|}{|
{\operatorname{det}}(\Xi_j)|^2}\|f\|^2_{\mathcal{N}_{\Phi
_{j-1}}(\mathbb{R}^d)}\nonumber
%
\end{eqnarray}
for $1\le j\le J$ where $c_1$ and $c_2$ satisfy (\ref{c1c2}), and
$\Phi
_{j-1}$ and $\Phi_j$ satisfy relations~(\ref{kernelRelation1}) and~(\ref
{kernelRelation2}).
\end{prop}
\begin{pf}
If $f\notin\mathcal{N}_{\Phi_{j-1}}(\mathbb{R}^d)$, then $\|f\|
^2_{\mathcal{N}_{\Phi_{j-1}}(\mathbb{R}^d)}=\infty$ and (\ref
{NormRelationEquation}) is true. Now, assume $f\in\mathcal{N}_{\Phi
_{j-1}}(\mathbb{R}^d)$,
and note that
\[
\frac{\omega'\Xi'_j\Xi_j\omega}{\omega'\Xi'_{j-1}\Xi
_{j-1}\omega}\le
\lambda_{\max}(\Theta'_{j-1}\Theta_{j-1}\Xi'_j\Xi_j).
\]
If $\lambda_{\max}(\Theta'_{j-1}\Theta_{j-1}\Xi'_j\Xi_j)\le1$, then
%
\begin{eqnarray}\label{FourierInequality}
&&\omega'\Xi'_j\Xi_j\omega\le\omega'\Xi'_{j-1}\Xi_{j-1}\omega
\nonumber\\
&&\quad\Longrightarrow\quad\frac{1}{c_1}\hat{\Phi}(\Xi_j\omega)\ge\hat
{\Upsilon
}(\Xi_j\omega)\ge\hat{\Upsilon}(\Xi_{j-1}\omega)\ge\frac
{1}{c_2}\hat
{\Phi}(\Xi_{j-1}\omega)\nonumber\\[-8pt]\\[-8pt]
&&\quad\Longrightarrow\quad\frac{1}{\hat{\Phi}(\Xi_j\omega)^{2^{J-j}}}\le
\biggl(\frac{c_2}{c_1}\biggr)^{2^{J-j}}\frac{1}{\hat{\Phi}(\Xi
_{j-1}\omega
)^{2^{J-j}}}\nonumber\\
&&\quad\Longrightarrow\quad\frac{1}{\hat{\Psi}^{J-j}(\Xi_{j}\omega)}\le
\biggl(\frac
{c_2}{c_1}\biggr)^{2^{J-j}}\frac{1}{\hat{\Psi}^{J-j}(\Xi
_{j-1}\omega)},
\nonumber
\end{eqnarray}
where the first implication follows from (\ref{c1c2}), the second
implication follows since the right- and left-hand sides are positive
and the final implication follows from the relations (\ref
{kernelRelation1}) and (\ref{kernelRelation2}) and the properties of
Fourier transforms of convolutions. So,
\begin{eqnarray*}
&&\|f\|^2_{\mathcal{N}_{\Phi_{j-1}}(\mathbb{R}^d)}\\
&&\qquad=\|f\|^2_{\mathcal
{N}_{\Psi^{J-(j-1)}_{\Theta_{j-1}}}(\mathbb{R}^d)}\\
&&\qquad=(2\pi)^{-d/2}\int_{\mathbb{R}^d}\frac{|\hat{f}(\omega)|^2}{\hat
{\Psi
}^{J-(j-1)}_{\Theta_{j-1}}(\omega)}\,\mathrm{d}\omega\\
&&\qquad=\frac{(2\pi)^{-d/2}}{|{\operatorname{det}}(\Xi_{j-1})|}\int_{\mathbb
{R}^d}\frac
{|\hat{f}(\omega)|^2}{\hat{\Psi}^{J-(j-1)}(\Xi_{j-1}\omega)}
\,\mathrm{d}\omega\\
&&\qquad=\frac{(2\pi)^{-d/2}}{|{\operatorname{det}}(\Xi_{j-1})|}\int_{\mathbb
{R}^d}\frac
{|\hat{f}(\omega)|^2}{\widehat{\Psi^{J-j}\ast\Psi^{J-j}}(\Xi
_{j-1}\omega
)}\,\mathrm{d}\omega\\
&&\qquad=\frac{(2\pi)^{-d}}{|{\operatorname{det}}(\Xi_{j-1})|}\int_{\mathbb
{R}^d}\frac
{|\hat{f}(\omega)|^2}{\hat{\Psi}^{J-j}(\Xi_{j-1}\omega)^2}
\,\mathrm{d}\omega
\\
&&\qquad=(2\pi)^{-d}\frac{|{\operatorname{det}}(\Xi_j)|^2}{|{\operatorname{det}}(\Xi
_{j-1})|}\int
_{\mathbb{R}^d}\frac{|\hat{f}(\omega)|^2}{|{\operatorname{det}}(\Xi_j)|^2\hat
{\Psi
}^{J-j}(\Xi_{j-1}\omega)^2}\,\mathrm{d}\omega\\
&&\qquad\ge(2\pi)^{-d}\frac{|{\operatorname{det}}(\Xi_j)|^2}{|{\operatorname{det}}(\Xi
_{j-1})|}
\biggl(\frac{c_1}{c_2}\biggr)^{2^{J-j+1}}\int_{\mathbb{R}^d}\frac{|\hat
{f}(\omega)|^2}{|{\operatorname{det}}(\Xi_j)|^2\hat{\Psi}^{J-j}(\Xi_{j}\omega
)^2}\,\mathrm{d}\omega\\
&&\qquad=(2\pi)^{-d/2}\frac{|{\operatorname{det}}(\Xi_j)|^2}{|{\operatorname{det}}(\Xi
_{j-1})|}
\biggl(\frac{c_1}{c_2}\biggr)^{2^{J-j+1}}\int_{\mathbb{R}^d}\frac{|\hat
{f}(\omega)|^2}{\widehat{\Psi^{J-j}_{\Theta_j}\ast\Psi
^{J-j}_{\Theta
_j}}(\omega)}\,\mathrm{d}\omega\\
&&\qquad=\frac{|{\operatorname{det}}(\Xi_j)|^2}{|{\operatorname{det}}(\Xi_{j-1})|}\biggl(\frac
{c_1}{c_2}\biggr)^{2^{J-j+1}}\|f\|^2_{\mathcal{N}_{\Phi_j\ast\Phi
_j}(\mathbb{R}^d)},
\end{eqnarray*}
where the first equality follows from relation (\ref{kernelRelation2}),
the second equality follows from the inner product representation (\ref
{InnerProdRd}), the third equality follows from the scaled Fourier
transform relation (\ref{FourierRelationship}), the fourth equality
follows from the definition of $\Psi^{J-(j-1)}$ (\ref
{kernelRelation1}), the fifth equality follows from the properties of
Fourier transforms of convolutions, the sixth equality follows by
multiplying by $|{\operatorname{det}}(\Xi_j)|^2/|{\operatorname{det}}(\Xi_j)|^2$, the
inequality follows from the development~(\ref{FourierInequality}), the
seventh equality follows from the scaled Fourier transform relation~(\ref{FourierRelationship}) and the properties of Fourier transforms of
convolutions and the final equality follows from the inner product
representation (\ref{InnerProdRd}).
\end{pf}
%

In most applications, the domain of interest $\Omega$ is a strict
subset of $\mathbb{R}^d$. If
$f\in\mathcal{N}_{\Phi_{1}}(\Omega)$, then $f$ can be extended to
$Ef\in
\mathcal{N}_{\Phi_{1}}(\mathbb{R}^d)$
\cite{wendland2005}
with
%
\begin{eqnarray}\label{extension}
\|f\|_{\mathcal{N}_{\Phi_{1}}(\Omega)}&=&\|Ef\|_{\mathcal{N}_{\Phi
_{1}}(\mathbb{R}^d)},\nonumber\\[-8pt]\\[-8pt]
\|f\|_{\mathcal{N}_{\Phi_{2}}(\Omega)}&\le&\|Ef\|_{\mathcal
{N}_{\Phi
_{2}}(\mathbb{R}^d)}\nonumber
\end{eqnarray}
for all $\Phi_2$. Combining (\ref{extension}) with Proposition
\ref{NormRelationProposition} gives the following corollary.

\begin{corollary}\label{NormRelationPropositionOmega}
If the assumptions of Proposition \ref{NormRelationProposition} are
satisfied, then
%
\begin{eqnarray}\label{NormRelationEquationOmega}
&&\lambda_{\max}(\Theta'_{j-1}\Theta_{j-1}\Xi'_j\Xi_j)\le1\nonumber\\[-8pt]\\[-8pt]
&&\quad\Longrightarrow\quad\|f\|^2_{\mathcal{N}_{\Phi_j\ast\Phi_j}(\Omega
)}\le
\biggl(\frac{c_2}{c_1}\biggr)^{2^{J-(j-1)}}\frac{|{\operatorname{det}}(\Xi
_{j-1})|}{|{\operatorname{det}}(\Xi_j)|^2}\|f\|^2_{\mathcal{N}_{\Phi
_{j-1}}(\Omega)}.
\nonumber
%
\end{eqnarray}
\end{corollary}
\begin{pf}
If $f\notin\mathcal{N}_{\Phi_{j-1}}(\Omega)$, then $\|f\|
^2_{\mathcal
{N}_{\Phi_{j-1}}(\Omega)}=\infty$ and (\ref{NormRelationEquationOmega})
is true. Now, assume $f\in\mathcal{N}_{\Phi_{j-1}}(\Omega)$ and extend
$f$ to $Ef\in\mathcal{N}_{\Phi_{j-1}}(\mathbb{R}^d)$ with $\|Ef\|
^2_{\mathcal{N}_{\Phi_{j-1}}(\mathbb{R}^d)}=\|f\|^2_{\mathcal
{N}_{\Phi
_{j-1}}(\Omega)}$. Then
\begin{eqnarray*}
\|f\|^2_{\mathcal{N}_{\Phi_j\ast\Phi_j}(\Omega)}&\le&\|Ef\|
^2_{\mathcal
{N}_{\Phi_j\ast\Phi_j}(\mathbb{R}^d)}\\
&\le&\biggl(\frac{c_2}{c_1}\biggr)^{2^{J-(j-1)}}\frac{|{\operatorname{det}}(\Xi
_{j-1})|}{|{\operatorname{det}}(\Xi_j)|^2}\|Ef\|^2_{\mathcal{N}_{\Phi
_{j-1}}(\mathbb{R}^d)}\\
&=&\biggl(\frac{c_2}{c_1}\biggr)^{2^{J-(j-1)}}\frac{|{\operatorname{det}}(\Xi
_{j-1})|}{|{\operatorname{det}}(\Xi_j)|^2}\|f\|^2_{\mathcal{N}_{\Phi
_{j-1}}(\Omega
)},
%
\end{eqnarray*}
where the first inequality follows from (\ref{extension}), the second
inequality follows from Proposition \ref{NormRelationProposition} and
the equality follows from the property of the chosen extension.
\end{pf}
%

\subsection{Error bound for multi-step interpolator}

Combining Theorem \ref{nativeSpaceBound} with Corollary \ref
{NormRelationPropositionOmega}, we are able to obtain the following
theorem bounding the native space norm of the multi-step interpolator's error.
%
\begin{theorem}\label{nativeSpaceErrorBound}
Under the assumptions of Theorem \ref{nominalResult} and Proposition~\ref{NormRelationProposition},
\[
\Biggl\|f-\sum_{j=1}^J\mathcal{P}^j\Biggr\|_{\mathcal{N}_{\Phi_J}(\Omega)}\le
C_{\Phi,J}\|f\|_{\mathcal{N}_{\Phi_0}(\Omega)}\prod_{j=1}^J\biggl\{
\frac
{\sqrt{|{\operatorname{det}}(\Xi_{j-1})|}}{|{\operatorname{det}}(\Xi_j)|}(\|\Theta
_j\|
^k_2 h^k_{X_j})^{2^{J-j-1}}\biggr\}.
\]
\end{theorem}
\begin{pf}
First applying Theorem \ref{nativeSpaceBound} and then applying
Proposition \ref{NormRelationProposition} gives
\begin{eqnarray*}
\Biggl\|f-\sum_{j=1}^J\mathcal{P}^j\Biggr\|_{\mathcal{N}_{\Phi_J}(\Omega)}&\le&
C_\Phi\|\Theta_J\|_2^{k/2}h_{X_J}^{k/2}\Biggl\|f-\sum_{j=1}^{J-1}\mathcal
{P}^j\Biggr\|_{\mathcal{N}_{\Phi_J\ast\Phi_J}(\Omega)}\\
&\le& C_{\Phi,J}\|\Theta_J\|_2^{k/2}h_{X_J}^{k/2}\frac{\sqrt{|
{\operatorname{det}}(\Xi_{J-1})|}}{|{\operatorname{det}}(\Xi_J)|}\Biggl\|f-\sum_{j=1}^{J-1}\mathcal
{P}^j\Biggr\|_{\mathcal{N}_{\Phi_{J-1}}(\Omega)}.
\end{eqnarray*}
For $J\ge2$, repeat the above argument $J-1$ more times, and note that
$\Phi_{J-j}$ has $k2^{j}$ continuous derivatives.
\end{pf}

By applying Theorem \ref{nominalResult} to the error $f-\sum
_{j=1}^J\mathcal{P}^j$, an additional multiple of $h_{X_J}^{k/2}$ is
obtained in the following theorem.
%
\begin{theorem}\label{nominalStar}
Under the assumptions of Theorem \ref{nominalResult} and Proposition~\ref{NormRelationProposition},
\begin{eqnarray*}
&&\Biggl|f(x)-\sum_{j=1}^J\mathcal{P}^j(x)\Biggr|\nonumber\\[-8pt]\\[-8pt]
&&\qquad\le C_{\Phi,J}\|f\|_{\mathcal
{N}_{\Phi_0}(\Omega)}\|\Theta_J\|_2^{k/2}h_{X_J}^{k/2}\prod
_{j=1}^J
\biggl\{\frac{\sqrt{|{\operatorname{det}}(\Xi_{j-1})|}}{|{\operatorname{det}}(\Xi_j)|}(\|
\Theta
_j\|^k_2 h^k_{X_j})^{2^{J-j-1}}\biggr\}.
\end{eqnarray*}
\end{theorem}

\section{Examples}

\begin{figure}

\includegraphics{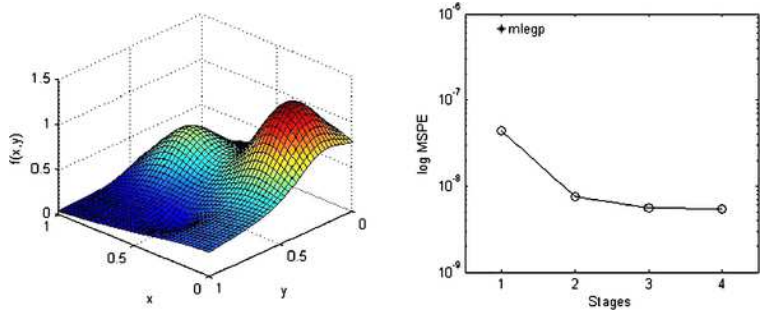}

\caption{Left panel: Franke's function. Right panel:
log mean squared prediction error versus number of stages (circles) and
using \texttt{mlegp} (asterisk).}\label{multiscaleexample}
\end{figure}

First, consider using the multi-step procedure to interpolate Franke's function
\begin{eqnarray*}
f(x,y)&=&\tfrac{3}{4}\operatorname{exp}\bigl\{-\bigl((9x-2)^2+(9y-2)^2\bigr)/4\bigr\}\\
&&{} +\tfrac{3}{4}\operatorname{exp}\bigl\{-\bigl((9x+1)^2/49-(9y+1)^2/10\bigr)\bigr\}\\
&&{} +\tfrac{1}{2}\operatorname{exp}\bigl\{-\bigl((9x-7)^2+(9y-3)^2\bigr)/4\bigr\}\\
&&{} -\tfrac{1}{5}\operatorname{exp}\bigl\{-\bigl((9x-4)^2+(9y-7)^2\bigr)\bigr\}
\end{eqnarray*}
%
shown in the left panel of Figure \ref{multiscaleexample}.
Theorems \ref{numericThm} and \ref{thmg} indicate that each of the
nested data sets should have well-separated points in the full
%
%
dimension as well as lower-dimensional projections to give small
numeric error, and Theorem \ref{nominalStar} indicates that each of the
nested data sets should have small data-free regions in the full
dimension as well as lower-dimensional projections to give small
nominal error.
Training data are collected from Franke's function using a randomized
$(0,4,2)$-net in base 5 \cite{owen1995} with $5^4=625$ points, which
has a convenient nested structure with both the full and each
sub-design having small data-free regions and relatively well-spread
points in both the full and projected space, making it ideal for the multi-step
procedure.
Theorem \ref{thmg} indicates that a less smooth underlying kernel
$\Phi
$ will give more numerically accurate results, while Theorem \ref
{nominalStar} indicates that a~more smooth kernel will give more
nominally accurate results.
To balance these opposing forces in this moderately sized example, the
selected $\Phi$ is Wendland's
compactly supported kernel with four continuous derivatives~\cite
{wendland2005},
\begin{eqnarray*}
\Phi(x-y)&=&\phi\bigl(\sqrt{(x-y)'(x-y)}\bigr),\\
\phi(r)&=&(1-r)^{l+2}_+[(l^2+4l+3)r^2+(3l+6)r+3],\qquad
l=\lfloor
d/2\rfloor+3,
\end{eqnarray*}
and the rescaling matrices $\Theta_1,\ldots,\Theta_J$
are restricted to be diagonal, so each input is re-scaled separately.
The re-scalings for each stage
are chosen by leave-one-out cross-validation, for which a simple
short-cut formula holds making
computation undemanding for this moderately sized problem, although
$A^{-1}_{X_j,\Phi_j}$ needs to\vspace*{1pt} be calculated.
In particular, the $i$th cross-validation error at stage $j$ is
\cite{rippa}
%
\begin{equation}\label{shortcut}
e_{(i)}=\frac{\alpha^j_i}{B^j_{ii}}, \qquad B^j=A^{-1}_{X_j,\Phi_j}.
\end{equation}
In this example, the single-stage sample size is $n_1=625$, the
two-stage sample sizes are $n_1=250$ and $n_2=625$, the three-stage
sample sizes are $n_1=250$, $n_2=375$ and $n_3=625$ and the four-stage
sample sizes are $n_1=250$, $n_2=375$, $n_3=500$ and $n_4=625$.
The nested data sets are $X_j=\{x_i\in X\dvtx i\le n_j\}$.
The right panel of Figure \ref{multiscaleexample} shows the logarithm
of the mean squared prediction error on a test set of 1,000 randomly
generated uniform points.
Notice\vspace*{1pt} that the mean squared prediction error is improved from
$4.4\times10^{-8}$ to $5.4\times10^{-9}$.
A Gaussian process fit using the \texttt{mlegp} package \cite{mlegp} in
\texttt{R}, on the other hand, has mean squared prediction error
$6.8\times10^{-7}$.

\begin{figure}

\includegraphics{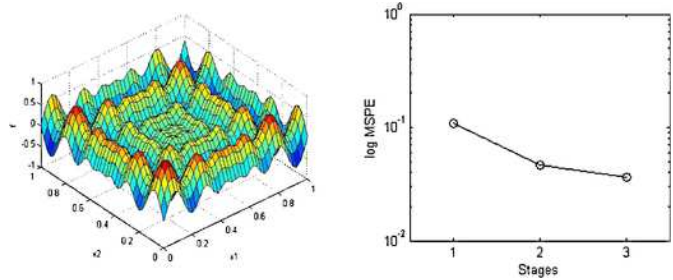}

\caption{Left panel: two-dimensional projection of Schwefel's
function. Right panel: log mean squared prediction error versus
number of stages.}\label{multiscaleexample2}
\end{figure}

Next, consider using the multi-step procedure to interpolate Schwefel's
function for $d=5$
\[
f(x)=-\sum_{j=1}^d (1\mbox{,}000x_j-500)\sin\bigl(\sqrt
{|1\mbox{,}000x_j-500|}\bigr)/1\mbox{,}000,
\]
a two-dimensional projection of which with the remaining variables fixed
at~$1/2$ is shown in the left panel of Figure \ref{multiscaleexample2}.
This function is relatively complex and a very large training set is
needed to build an accurate emulator.
%
%
To ensure easy nesting and good space-filling properties for sub-designs,
data are collected from Schwefel's function using a randomized
$(0,8,5)$-net in base~5 with $5^8=390\mbox{,}625$ points.
In this example there is a great deal of potential for numeric problems so
Wendland's continuous,
compactly supported kernel,
\begin{eqnarray*}
\Phi(x-y)&=&\phi\bigl(\sqrt{(x-y)'(x-y)}\bigr),\\
\phi(r)&=&(1-r)^{l+2}_+,\qquad l=\lfloor d/2\rfloor+1
\end{eqnarray*}
with relatively little smoothness is selected.
The re-scaling matrices $\Theta_1,\ldots,\allowbreak\Theta_J$
are chosen to be \textit{fixed} scalar multiples of the identity,
$\Theta
_j=\theta_j I_d$, with
%
\begin{equation}\label{fixedTheta}
\theta_j=\biggl(\frac{n_j^2\pi^{d/2}}{10^7\Gamma(d/2+1)}\biggr)^{1/d},
\end{equation}
which ensures that each intepolation matrix $A_{X_j,\Phi_j}$ has less
than $10^7$ nonzero entries.
Edge effects in the five-dimensional cube ensure that the number of
nonzero entries is substantially less than $10^7$.
In this example, the single-stage sample size is $n_1=390\mbox{,}625$, the
two-stage sample sizes are $n_1=5^7=78\mbox{,}125$ and $n_2=390\mbox{,}625$ and the
three-stage sample sizes are $n_1=78\mbox{,}125$, $n_2=2\times5^7=156\mbox{,}250$ and
$n_3=390\mbox{,}625$.
The nested data sets are $X_j=\{x_i\in X\dvtx i\le n_j\}$.
The right panel of Figure \ref{multiscaleexample2} shows the logarithm
of the mean squared prediction error on a test set of 10,000 randomly
generated uniform points.
Notice that the mean squared prediction error is improved from $0.11$
to $0.036$.
On the other hand, the \texttt{mlegp} package runs out of memory trying
to fit a GP.

\section{Discussion}

We have presented the intuitively appealing and practically useful
multi-step interpolation procedure. This procedure is easy to use and
offers substantial improvements in overall accuracy in the emulation of
large-scale computer experiments. We introduced a decomposition of the
error of \textit{any} interpolator into nominal and numeric portions.
This decomposition is important because it allows the two sources of
error to be analyzed separately while emphasizing the interplay between
the two types of errors. We proved a very general result bounding the
numeric error of a~multi-step interpolator, of which an ordinary
interpolator is a special case. This result constitutes the only
complete and rigorous bound on the numeric error of the multi-step
interpolator. We proved that in the situation where the earlier stage
kernels are convolutions of the later stage kernels, substantial
nominal improvements can be realized. In the context of the multi-step
interpolator, this result is the most general and explicit of its kind.

Further work on the multi-step interpolation method will be explored in
the following directions.
First, its implementation details, along with
various examples, will be reported in a subsequent article,
to illustrate the theoretical results derived here.
The implementation of the method
requires the generation of nested data sites, for which the typical
choice in applied mathematics is nested grids.
Nested space-filling designs \cite{qian20091,qa10,qian20092,haaland1},
originally constructed for running multiple computer experiments with
different levels of accuracy,
are a better choice because of their good uniformity properties. Such
designs can be generated
by exploiting \textit{nesting} in orthogonal arrays~\cite
{HedayatSloaneStufken1999}, U designs \cite{tang93,tang94},
orthogonal Latin hypercubes \mbox{\cite{ye1998,sl06,bst09,lin09,lin10}} or
scrambled nets~\cite{owen1995}.
Second, emulation of computer models with qualitative and quantitative
factors is currently
getting increasing attention \cite{qww08,qian20094,han}. 
We plan to extend the multi-step procedure to accommodate these two
types of factors.
Third, beyond emulation of computer experiments,
singularity issues arise in fitting many other large kernel models.
We plan to introduce a~general multi-step framework for
fitting kernel based classification and regression methods with a large
number of observations.
As in the multi-step interpolation procedure,
this framework obtains nested data sites and then fits a kernel model
in multiple
steps, where in each step interpolation is
replaced by an appropriate procedure for the given problem.
New theoretical bounds on the nominal and numeric accuracy,
analogous to those in Sections~\ref{secnum} and~\ref{nominalsection},
will be derived for this framework. The required well-spread nested
data sites for the framework will be generated by using
nested space-filling designs or the efficient \textit{thinning} algorithm
\cite{floater} for observational data.
In the revision of this paper, we became aware of new theoretical developments
of the multi-step method in applied mathematics, including \cite{GiSlWe10}
and \cite{Wen10}.

\section*{Acknowledgments}

The authors thank the Editor, an Associate Editor and two referees for
their comments that have led to improvements in the paper. They also
thank Greg Fasshauer, Yizhi Zhang, Fred Hickernell and Grace Wahba for
their comments.


%

\printaddresses

\end{document}